\documentclass[11 pt]{amsart}
\usepackage{latexsym,amscd,amssymb, graphicx, amsthm, young}  
\usepackage{blkarray}
\usepackage{tikz-cd}

\usepackage[margin=1in]{geometry}

\newtheorem{theorem}{Theorem}

\newtheorem{conjecture}[theorem]{Conjecture}

\theoremstyle{definition}

\newcommand{\Frob}{{\mathrm {Frob}}}

\newcommand{\symm}{{\mathfrak{S}}}

\newcommand{\CC}{{\mathbb {C}}}

\newcommand{\LIS}{\mathrm{LIS}}

\renewcommand{\P}{\mathbb{P}}

\begin{document}

\title[Increasing Subsequences]
{Increasing Subsequences and Kronecker Coefficients}

\author{Jonathan Novak}
\author{Brendon Rhoades}
\address
{Department of Mathematics \newline \indent
University of California, San Diego \newline \indent
La Jolla, CA, 92093-0112, USA}
\email{(jinovak, bprhoades)@ucsd.edu}

\maketitle

A permutation $\pi$ in the symmetric group $\symm_n$ is said to have an increasing
subsequence of length $k$ if there are numbers $1 \leq i_1< \dots < i_k \leq n$ such that

	\begin{equation*}
		\pi(i_1) < \dots < \pi(i_k).
	\end{equation*}
	
\noindent
Similarly, $\pi$ has a decreasing subsequence of length $l$ if we can find $1 \leq j_1< \dots < j_l \leq n$ with

	\begin{equation*}
		\pi(j_1) > \dots > \pi(j_l).
	\end{equation*}
	
\noindent
Increasing and decreasing subsequences in permutations have been the subject of sustained interest in combinatorics since
1935, when Erd\H{o}s and Szekeres proved that, given any $k,l \in \mathbb{N},$ every permutation in $\mathfrak{S}_n$
contains either an increasing subsequence of length $k$ or a decreasing subsequence of length $l$ as soon as $n$ exceeds $(k-1)(l-1).$
This is a permutation version of Ramsey's theorem in which the threshold between possible disorder and certain order is 
explicitly computable.

In view of the Erd\H{o}s-Szekeres theorem, it is natural to wonder about the typical length of longest 
monotone subsequences in permutations. Since reversing a permutation interchanges its increasing
and decreasing subsequences, we may focus on just one of the two orientations. 
Let $\LIS_n$ denote the length of the longest increasing 
subsequence in a uniformly random sample from 
$\symm_n.$ What can be said about the distribution of $\LIS_n$? 
This question has been answered in the $n \to \infty$ limit: we have both a Law of Large Numbers
and a Central Limit Theorem for $\LIS_n.$ The LLN was was obtained by Vershik and Kerov in 1977, who showed
that 

	\begin{equation*}
		\lim_{n \to \infty} \frac{\LIS_n}{\sqrt{n}} = 2,
	\end{equation*}
	
\noindent
the convergence being in probability. The corresponding CLT is a breakthrough 1999
result of Baik, Deift, and Johansson, who discovered that

	\begin{equation*}
		\lim_{n \to \infty} \P \left( \frac{\LIS_n - 2\sqrt{n}}{n^{1/6}} \leq t \right) = F(t),
	\end{equation*}
	
\noindent
where $F(t)$ is the Tracy-Widom distribution from random matrix theory. For an informative discussion of these 
results, together with pertinent references, we refer the reader to Stanley's 2006 ICM contribution \cite{Stanley:ICM}.

The distribution of $\LIS_n$ at finite $n$ remains poorly understood --- even the question of unimodality has 
not been settled. Given that the distribution of most natural combinatorial statistics is unimodal, it is not unreasonable to 
conjecture that this holds true for $\LIS_n.$ The above limit theorems certainly suggest that this is the case, since they
tell us that when $n$ is large the distribution of $\LIS_n$ is concentrated in a small neighborhood of $2\sqrt{n}.$ 
However, this does not rule out the possibility of irregular ``noisy'' behavior before the limit. 
In combinatorics and elsewhere, the source of unimodality is often a stronger property: log-concavity \cite{Stanley:NY}.
Conjecturally, this is the case for $\LIS_n.$

\begin{conjecture}[Chen \cite{Chen}]
\label{ChenConjecture}
	For any $n \geq 3,$ the distribution of $\LIS_n$ is log-concave. That is, we have
	
		\begin{equation}
		\label{ineq:LogConcave}
			a_{n,k-1} a_{n,k+1} \leq a_{n,k}^2,
		\end{equation}
		
	\noindent
	for all $2 \leq k \leq n-1,$ where $a_{n,k}$ is the number of permutations in $\symm_n$ whose
	longest increasing subsequence has length $k.$
\end{conjecture}
	
The inequality \eqref{ineq:LogConcave} is equivalent to the existence of an injection 

	\begin{equation}
	\label{ineq:LogConcaveSets}
		\mathfrak{S}_{n,k-1} \times \mathfrak{S}_{n,k+1} \longrightarrow \mathfrak{S}_{n,k} \times \mathfrak{S}_{n,k},
	\end{equation}
	
\noindent 
where $\mathfrak{S}_{n,k}$ is the set of permutations in $\symm_n$ whose longest increasing 
subsequence has length equal to $k.$ 
It seems difficult to construct such an injection in any uniform way: tampering with a permutation
alters its increasing subsequences in complicated ways which are difficult to track. Instead, one 
may invoke the Robinson-Schensted correspondence \cite{Schensted}, which tells us that $a_{n,k}$ is equal to 
the number of pairs of standard Young tableaux on a common shape $\lambda \vdash n$ 
satisfying $\ell(\lambda)=k.$ Thus \eqref{ineq:LogConcaveSets} is equivalent to the existence of an
injection

	\begin{equation}
	\label{ineq:LogConcavePairs}
		\mathfrak{P}_{n,k-1} \times \mathfrak{P}_{n,k+1} \longrightarrow \mathfrak{P}_{n,k} \times \mathfrak{P}_{n,k},
	\end{equation}
	
\noindent
where $\mathfrak{P}_{n,k}$ is the set of Robinson-Schensted pairs with $n$ cells and $k$ rows.
Since operations on Young tableaux are easier to visualize than operations on permutations, such an injection may be easier to construct.
Another possibility would be to use the hook-length formula $f^\lambda=n!/H_\lambda$ for the number $f^\lambda$ of SYT of shape $\lambda,$ 
which combined with \eqref{ineq:LogConcavePairs} gives

	\begin{equation}
	\label{ineq:LogConcaveHooks}
		\sum_{\substack{\lambda \vdash n \\ \ell(\lambda)=k-1}}
		\sum_{\substack{\mu \vdash n \\ \ell(\mu)=k+1}} \frac{1}{H_\lambda^2 H_\mu^2} \leq 
		\sum_{\substack{\lambda \vdash n \\ \ell(\lambda)=k}}
		\sum_{\substack{\mu \vdash n \\ \ell(\mu)=k}} \frac{1}{H_\lambda^2 H_\mu^2}.
	\end{equation}
	
\noindent
This inequality is equivalent to \eqref{ineq:LogConcave}, but may be easier to work with. So far,
neither of these strategies has been successfully implemented; see \cite{BLS} for some partial results.

When faced with an intractable problem, there is nothing to lose and everything to gain by trying to solve
an even harder problem. We submit that the ``right'' way to approach Conjecture \ref{ChenConjecture} is to focus
on the following stronger conjecture. Let $R(\symm_n)$ be the representation ring of $\symm_n,$ i.e. the
commutative ring generated by the isomorphism classes $V^\lambda$ of irreducible complex representations of $\symm_n$ 
with operations being direct sum and tensor product of $\symm_n$-modules. Let us replace the sets 
$\symm_{n,k}$ appearing in \eqref{ineq:LogConcaveSets} with the modules $V_{n,k} \in R(\symm_n)$ defined by

	\begin{equation*}
		V_{n,k} = \bigoplus_{\substack{\lambda \vdash n \\ \ell(\lambda) = k}} f^\lambda V^\lambda.
	\end{equation*}
	
\noindent
Since $\dim V^\lambda=f^\lambda,$ we have

	\begin{equation*}
		\dim V_{n,k} = \sum_{\substack{\lambda \vdash n \\ \ell(\lambda) = k}} (f^\lambda)^2 = a_{n,k},
	\end{equation*}
	
\noindent
whence inequality \eqref{ineq:LogConcave} is equivalent to 

	\begin{equation}
	\label{eqn:LogConcaveDimensions}
		\dim V_{n,k-1} \dim V_{n,k+1} \leq \dim V_{n,k} \dim V_{n,k}.
	\end{equation}
	
\noindent
We propose the following strengthening of Conjecture \ref{ChenConjecture}.

	\begin{conjecture}
	\label{NovakRhoadesEquivariant}
	For any $n \geq 3,$ there exists an $\symm_n$-equivariant injection
	
	\begin{equation}
	\label{ineq:LogConcaveModules}
		V_{n,k-1} \otimes V_{n,k+1} \longrightarrow V_{n,k} \otimes V_{n,k}.
	\end{equation}
	
	\noindent
	for all $2 \leq k \leq n-1.$
	\end{conjecture}
	
Conjecture \ref{NovakRhoadesEquivariant} may equivalently be formulated as a numerical refinement of 
\eqref{ineq:LogConcaveHooks}. To see this, we make use of the Kronecker coefficients, which linearize multiplication 
in $R(\symm_n)$:	

	\begin{equation*}
		V^\lambda \otimes V^\mu = \bigoplus_{\nu \vdash n} g_{\lambda\mu}^\nu V^\nu.
	\end{equation*}

\noindent
Decomposing the source and target in \eqref{ineq:LogConcaveModules} into irreducibles, 
Conjecture \ref{NovakRhoadesEquivariant} claims the existence of an injective $\symm_n$-module
homomorphism

	\begin{equation*}
		\bigoplus_{\nu \vdash n} \left( \sum_{\substack{\lambda \vdash n \\ \ell(\lambda)=k-1}}
		\sum_{\substack{\mu \vdash n \\ \ell(\mu)=k+1}} f^\lambda f^\mu g_{\lambda\mu}^\nu \right) V^\nu
		\longrightarrow
		\bigoplus_{\nu \vdash n} \left( \sum_{\substack{\lambda \vdash n \\ \ell(\lambda)=k}}
		\sum_{\substack{\mu \vdash n \\ \ell(\mu)=k}} f^\lambda f^\mu g_{\lambda\mu}^\nu \right) V^\nu.
	\end{equation*}

\noindent
Thus, by Schur's lemma and the hook-length formula, Conjecture \ref{NovakRhoadesEquivariant} is equivalent
to the following numerical inequality.

	\begin{conjecture}
	\label{NovakRhoadesNumerical}
	For any $n \geq 3$ and $\nu \vdash n,$ we have

	\begin{equation}
	\label{ineq:LogConcaveKronecker}
	 \sum_{\substack{\lambda \vdash n \\ \ell(\lambda)=k-1}}
		\sum_{\substack{\mu \vdash n \\ \ell(\mu)=k+1}} \frac{g_{\lambda\mu}^\nu}{H_\lambda H_\mu}
		\leq \sum_{\substack{\lambda \vdash n \\ \ell(\lambda)=k}}
		\sum_{\substack{\mu \vdash n \\ \ell(\mu)=k}} \frac{g_{\lambda\mu}^\nu}{H_\lambda H_\mu}
	\end{equation}
	
	\noindent
	for all $2 \leq k \leq n-1.$
	\end{conjecture}
	
\noindent
The inequality \eqref{ineq:LogConcaveHooks} is recovered from Conjecture \ref{NovakRhoadesNumerical} by summing
\eqref{ineq:LogConcaveKronecker} over all $\nu \vdash n.$
	
Yet another equivalent formulation of Conjecture \ref{NovakRhoadesEquivariant} may be obtained
by means of the Frobenius isomorphism 

	\begin{equation*}
		\Frob \colon R(\symm_n) \longrightarrow \Lambda_n,
	\end{equation*}
	
\noindent
where $\Lambda_n$ is the ring of homogeneous symmetric functions of degree $n$ equipped with 
the Kronecker product. We recall that the Kronecker product in $\Lambda_n$ is defined via bilinear
extension of the rule
 
	\begin{equation}
		s_{\lambda} * s_{\mu} := \sum_{\nu \vdash n} g_{\lambda\mu}^\nu  s_{\nu},
	\end{equation}
	
\noindent
where $\{s_\lambda \colon \lambda \vdash n\}$ is the Schur function basis. 
The Frobenius isomorphism is defined by $\Frob(V^\lambda)=s_\lambda.$
In particular, the Frobenius image of the increasing subsequence module 
$V_{n,k} \in R(\symm_n)$ appearing in Conjecture \ref{NovakRhoadesEquivariant} is the symmetric function
$S_{n,k} \in \Lambda_n$ given by

	\begin{equation*}
		S_{n,k} = \sum_{\substack{\lambda \vdash n \\ \ell(\lambda)=k}} f^\lambda s_\lambda.
	\end{equation*}
	
\noindent
Given symmetric functions $F,G \in \Lambda_n,$ let us write $F \leq G$ if the difference
$G-F$ is Schur positive, i.e. if the coefficients $c_\lambda$ defined by 

	\begin{equation*}
		G-F = \sum_{\lambda \vdash n} c_\lambda s_\lambda
	\end{equation*}
	
\noindent
are nonnegative. Conjecture \ref{NovakRhoadesEquivariant} may then be restated as follows.

	\begin{conjecture}
	\label{NovakRhoadesSymmetric}
		For any $n \geq 3,$ we have
		
			\begin{equation*}
				S_{n,k-1} * S_{n,k+1} \leq S_{n,k} * S_{n,k}
			\end{equation*}
			
		\noindent
		for all $2 \leq k \leq n-1.$
	\end{conjecture}
	
\noindent
Conjecture \ref{NovakRhoadesSymmetric} is a useful equivalent formulation of 
Conjecture \ref{NovakRhoadesEquivariant} in that Schur positivity is a well-developed 
topic in algebraic combinatorics. Conjecture \ref{NovakRhoadesSymmetric} has been
verified for $n \leq 15$ on a computer.
 
How could one go about proving Conjecture \ref{NovakRhoadesEquivariant}? Let us 
illustrate how a successful argument might look by outlining a representation-theoretic
proof of a much simpler proposition: the log-concavity of binomial coefficients.
From the hook-length formula, it is clear that the dimension of the irreducible representation of $\symm_n$ corresponding to the 
hook $\lambda = (n-k, 1^{k})$ is given by ${n-1 \choose k}$ for all $0 \leq k \leq n-1$:
\begin{equation}
\dim V^{(n-k,1^k)} = {n-1 \choose k}.
\end{equation}
The difference ${n - 1 \choose k}^2 - {n - 1 \choose k-1} \cdot {n-1 \choose k+1}$ may therefore be expressed as a difference of dimensions of Kronecker products:
\begin{multline}
{n - 1 \choose k}^2 - {n - 1 \choose k-1} \cdot {n-1 \choose k+1} = \\
\dim (V^{(n-k,1^k)} \otimes V^{(k+1,1^{n-k-1})}) - \dim ( V^{(n-k+1,1^{k-1})} \otimes V^{(k+2,1^{n-k-2})}).
\end{multline}
The sequence of binomial coefficients with upper index $n-1$ will thus be certified log-concave if we can exhibit an $\symm_n$-module
$V_{n,k}$ whose Frobenius image is 
\begin{equation}
\label{v-characterization}
\Frob(V_{n,k}) = s_{(n-k,1^k)} * s_{(k+1,1^{n-k-1})} - s_{(n-k+1,1^{k-1})} * s_{(k+2, 1^{n-k-2})}.
\end{equation}

The required module was found by
Kim and Rhoades \cite{KR}.  Let $\theta_1, \dots, \theta_n, \xi_1, \dots, \xi_n$ be a list of $2n$ anticommuting
variables and consider the exterior algebra
\begin{equation}
\wedge \{ \Theta_n, \Xi_n\} := \wedge \{ \theta_1, \dots, \theta_n, \xi_1, \dots, \xi_n \}
\end{equation}
over $\CC$ generated by these variables.  This is a $\CC$-vector space of dimension $2^{2n}$ which carries the
bigrading 
\begin{equation}
\wedge \{ \Theta_n, \Xi_n \} = \bigoplus_{i, j = 0}^n \wedge \{ \Theta_n, \Xi_n \}_{i,j}
\end{equation}
induced by considering the degree of the $\theta$-variables and $\xi$-variables separately.
In Physics, anticommuting variables are called ``fermionic,'' and the relation $\theta_i^2 = 0$
corresponds to the Pauli Exclusion Principle: no two fermions may occupy the same state at the same time.
Consider the diagonal action of $\symm_n$  on $\wedge \{ \Theta_n, \Xi_n \}$, viz.
\begin{equation}
w \cdot \theta_i := \theta_{w(i)} \quad \quad w \cdot \xi_i := \xi_{w(i)} \quad \quad w \in \symm_n, \, \, 1 \leq i \leq n
\end{equation}
and denote by $\langle \wedge \{\Theta_n, \Xi_n \}^{\symm_n}_+ \rangle \subseteq \wedge \{\Theta_n, \Xi_n \}$
the two-sided ideal generated by $\symm_n$-invariants with vanishing constant term.
The {\em fermionic diagonal coinvariant} ring is defined in \cite{KR} by 
\begin{equation}
FDR_n := \wedge \{ \Theta_n, \Xi_n \} / \langle \wedge \{\Theta_n, \Xi_n \}^{\symm_n}_+ \rangle. 
\end{equation}
This is a doubly graded $\symm_n$-module and an anticommutative version of the 
{\em diagonal coinvariant ring} \cite{Haiman}.

\begin{theorem} (Kim-R. \cite{KR})
\label{fermion-theorem}
The $(i,j)$-graded piece $(FDR_n)_{i,j}$ is zero unless $i + j < n$. When $i + j < n$ we have
\begin{equation}
\Frob( (FDR_n)_{i,j} ) = s_{(n-i,1^i)} * s_{(n-j,1^j)} - s_{(n-i+1,1^{i-1})} * s_{(n-j+1,1^{j-1})},
\end{equation}
where we interpret $s_{(n-i+1,1^{i-1})} * s_{(n-j+1,1^{j-1})} = 0$ if $i = 0$ or $j = 0$.
In particular, for $0 \leq k \leq n-1$ if we set $V_{n,k} := (FDR_n)_{k,n-k-1}$ then 
$V_{n,k}$ has Frobenius image given by Equation~\eqref{v-characterization}.
\end{theorem}

Theorem~\ref{fermion-theorem} implies that the symmetric function in 
Equation~\eqref{v-characterization} is Schur-positive, and taking vector space dimensions yields
the log-concavity of the sequence ${n-1 \choose 0}, {n-1 \choose 1}, \dots, {n - 1 \choose n-1}$.

\section*{Acknowledgements}

J. Novak was partially supported by NSF Grant DMS-1812288 and a Lattimer Fellowship.
B. Rhoades was partially supported by NSF Grant DMS-1500838 and 
DMS-1953781.

\end{document}